\documentclass[12pt]{article}

\usepackage{enumerate}
\usepackage{amsmath, amsfonts, amstext, amssymb, amsthm, amsbsy}
\usepackage[T1]{fontenc} 
\newtheorem{ex}{Example}[section]
\newtheorem{theo}{Theorem}[section]
\newtheorem{theo*}{Theorem*}[section]
\newtheorem{prop}{Proposition}[section]
\newtheorem{cor}{Corollary}[section]
\newtheorem{df}{Definition}[section]
\newtheorem{lem}{Lemma}[section]
\newtheorem{rem}{Remark}[section]
\newenvironment{preuve}[1]
{\noindent{\textit{Proof #1:}}}{\hfill\textbf{$\square$}\newline }

\newcommand{\cC}{{\cal C}}

\newcommand{\cF}{{\cal F}}

\newcommand{\cL}{{\cal L}}

\newcommand{\cN}{{\cal N}}

\newcommand{\cP}{{\cal P}}

\newcommand{\cS}{{\cal S}}

\newfont{\msbm}{msbm10 scaled\magstep1}
\newfont{\msbms}{msbm7 scaled\magstep1} 

\newcommand{\bbC}{\mbox{$\mbox{\msbm C}$}}

\newcommand{\bbE}{\mbox{$\mbox{\msbm E}$}}

\newcommand{\bbN}{\mbox{$\mbox{\msbm N}$}}

\newcommand{\bbP}{\mbox{$\mbox{\msbm P}$}}

\newcommand{\bbR}{\mbox{$\mbox{\msbm R}$}}

\newcommand{\bbZ}{\mbox{$\mbox{\msbm Z}$}}

\newcommand{\bbsC}{\mbox{$\mbox{\msbms C}$}}

\newcommand{\bbsR}{\mbox{$\mbox{\msbms R}$}}

\newcommand{\bbsZ}{\mbox{$\mbox{\msbms Z}$}}

\begin{document}

\title{Strong mixing properties of max-infinitely divisible random fields}
\author{Cl\'ement Dombry\footnote{Universit\'e de
Poitiers, Laboratoire de Mathématiques et Applications, UMR CNRS 7348, T\'el\'eport 2, BP 30179, F-86962 Futuroscope-Chasseneuil cedex,
France.  Email: Clement.Dombry@math.univ-poitiers.fr}\ \ and Fr\'ed\'eric Eyi-Minko\footnote{Universit\'e de
Poitiers, Laboratoire de Mathématiques et Applications, UMR CNRS 7348, T\'el\'eport 2, BP 30179, F-86962 Futuroscope-Chasseneuil cedex,
France.  Email: Frederic.Eyi.minko@math.univ-poitiers.fr}}
\date{}
\maketitle

\begin{abstract}
Let $\eta=(\eta(t))_{t\in T}$  be a sample continuous max-infinitely random field on a locally compact metric space $T$. 
For a closed subset $S\in T$, we note  $\eta_{S}$  the restriction of $\eta$ to $S$. We consider $\beta(S_1,S_2)$ the absolute regularity coefficient between $\eta_{S_1}$ and $\eta_{S_2}$, where  $S_1,S_2$ are two disjoint closed subsets of $T$. Our main result is a simple upper bound for $\beta(S_1,S_2)$ involving  the exponent measure $\mu$  of $\eta$: we prove that $\beta(S_1,S_2)\leq 2\int \bbP[\eta\not<_{S_1} f,\ \eta\not <_{S_2} f]\,\mu(df)$,
where  $f\not<_{S} g$ means that there exists $s\in S$ such that $f(s)\geq g(s)$.

If $\eta$ is a simple max-stable random field, the upper bound is related to the so-called extremal coefficients: for countable disjoint sets $S_1$ and $S_2$, we obtain $\beta(S_1,S_2)\leq 4\sum_{(s_1,s_2)\in S_1\times S_2}(2-\theta(s_1,s_2))$, where $\theta(s_1,s_2)$ is the pair extremal coefficient.

As an application, we show that these new estimates entail a central limit theorem for stationary max-infinitely divisible random fields on $\bbZ^d$. In the stationary max-stable case, we derive the asymptotic normality of three simple estimators of the pair extremal coefficient.
\end{abstract}
\\

\section{Introduction}
Max-stable random fields turn out to be fundamental models for spatial 
extremes since they arise as the the limit of rescaled
maxima. More precisely, consider the component-wise maxima
\[
\eta_n(t)=\max_{1\leq i\leq n} \xi_i(t),\quad t\in T,
\]
of independent copies $\xi_i$, $i\geq 1$, of a
random field $\xi=(\xi(t))_{t\in T}$. If the random field
$\eta_n=(\eta_n(t))_{t\in T}$ converges in distribution, as
$n\to\infty$, under suitable affine normalization, then its limit
$\eta=\{\eta(t)\}_{t\in T}$ is necessarily max-stable. 
Therefore, max-stable random fields play a central role in extreme value theory,
just like Gaussian random fields do in the classical statistical theory based on the central limit
Theorem.

Max-stable processes have been studied extensively in the last decades and
many of their properties are well-understood. For example, the structure of their
finite dimensional distributions is well known and insightful Poisson point process
or spectral representations are available. Also the theory has been extended to 
max-infinitely divisible (max-i.d.) processes. See for example the seminal works by Resnick \cite{R08}, 
de Haan \cite{dH78,dH84}, de Haan and Pickands \cite{dHP86}, Giné Hahn and Vatan \cite{GHV90}, Resnick and Roy \cite{RR91} and many others. More details and further references can be found in the manographs by Resnick \cite{R08} or de Haan and Fereira \cite{dHF06}.

The questions of mixing and ergodicity of max-stable random process indexed by $\bbR$ or $\bbZ$ have been addressed  recently.  First results by Weintraub \cite{W91} in the max-stable case  have  been completed by Stoev \cite{S10}, providing  necessary and sufficient conditions for mixing of max-stable process based on their spectral representations. More recently, Kabluchko and Schlather \cite{KS10} extend these results and obtain necessary and sufficient conditions for both mixing and ergodicity of max-i.d. random process. They define the dependence function of a stationary max-i.d. random process $\eta=(\eta(t))_{t\in\bbsZ}$ by
\[
\tau_a(h)=\log\frac{\bbP[\eta(0)\leq a, \eta(h)\leq a]}{\bbP[\eta(0)\leq a]\bbP[\eta(h)\leq a]},\quad a>\mathrm{essinf}\  \eta(0),\quad h\in\bbZ.
\]
Then, it holds with  $\ell=\mathrm{essinf}\  \eta(0)$:
\begin{itemize}
\item $\eta$ is mixing  if and only if for all $a>\ell$,  $\tau_a(n)\to 0$ as $n\to+\infty$;
\item $\eta$ is ergodic  if and only if for all $a>\ell$,  $n^{-1}\sum_{h=1}^n \tau_a(h)\to 0$ as $n\to+\infty$.
\end{itemize}
Ergodicity is strongly connected to the strong law of large numbers via the ergodic theorem. The above results find natural applications in statistics to obtain strong consistency of several natural estimators based on non-independent but ergodic observations.

Going a step further, we address in this paper the issue of estimating the strong mixing coefficients of max-i.d. random fields. In some sense, ergodicity and mixing state that the restrictions $\eta_{S_1}$ and $\eta_{S_2}$ to two subsets $S_1,S_2$ become almost independent when the distance between $S_1$ and $S_2$ goes to infinity. Strong mixing coefficients  make this statement quantitative: we  introduce two standard  mixing coefficients $\alpha(S_1,S_2)$ and $\beta(S_1,S_2)$ that measure how much $\eta_{S_1}$ and $\eta_{S_2}$ differ from independence. The rate of decay of those coefficients as the distance between $S_1$ and $S_2$ goes to infinity is a crucial point for the central limit theorem (see Appendix \ref{sec:TCL}). As an application, we consider the asymptotic normality of three simple estimators of  the extremal coefficients of a simple max-stable stationary random field on $\bbZ^d$.

Our approach differ from those  of Stoev \cite{S10} based on spectral representations and  of Schlather and Kabluchko \cite{KS10} based on exponent measures. It relies on the Poisson point process representation of max-i.d. random fields offered by Hahn, Giné and Vatan \cite{GHV90} (see also Appendix \ref{sec:A1}) and on the notions of extremal and subextremal points recently introduced by the authors \cite{DEM11b}. Palm theory for Poisson point process and Slyvniak's formula are also a key tool (see Appendix \ref{sec:slyvniak}).

The structure of the paper is the following: the framework and results are detailed in the next Section; Section 3 is devoted to the proofs and an Appendix gathers some more technical details.
 
\section{Framework and results}

Let  $(\Omega,\cF,\bbP)$ be a probability space and $T$ be a locally compact metric space. 
We note $\bbC(T)=\bbC(T, [0,+\infty))$ the space of nonnegative functions endowed with the topology of uniform convergence on compact sets and $\cC$ its  Borel $\sigma$-field. Let $\mu$ be a  locally finite Borel measure  on $\bbC_0(T)=\bbC(T)\setminus\{0\}$ satisfying
\begin{equation}\label{eq:condmu}
\mu\big[\{f\in\bbC_0(T);\ \sup_K f>\varepsilon\}\big]<\infty\quad \mbox{for all compact }K\subset T \mbox{ and }\varepsilon>0,
\end{equation}
and $\Phi$ a Poisson point process (P.P.P.) on $\bbC_0(T)$ with intensity $\mu$. More rigorously, we should  consider $\Phi$ as a random point measure rather  than as a random set of points, since there may be points with multiplicities. It is however standard to consider $\Phi$ as a random set of points with possible repetitions. 

We consider the random process 
\[
\eta(t)=\max\{\phi(t),\ \phi\in\Phi\},\quad t\in T,
\]
with the convention that the maximum of the empty set is equal to $0$. Condition \eqref{eq:condmu} ensures that 
the random process $\eta$ is  continuous on $T$ (see \cite{GHV90} and Appendix~\ref{sec:A1}). Another property  is worth noting: $\eta$ is max-infinitely divisible. This means that for all $n\geq 1$, there exist  independent and identically distributed continuous random fields $(\eta_{i,n})_{1\leq i\leq n}$ such that
\[
\eta\stackrel{\cL}=\vee_{i=1}^n \eta_{i,n},
\]
where $\vee$ stands for pointwise maximum and $\stackrel{\cL}=$ for equality in distribution. Note that the max-infinite divisibility of $\eta$ is a simple consequence of the superposition Theorem for Poisson point processes. Furthermore, for all $t\in T$, the essential infimum of the random variable $\eta(t)$ is equal to $0$.  As shown by Gin\'e, Hahn and Vatan \cite{GHV90}, up to simple transformations, essentially  all  max-i.d. continuous random process on $T$ can be obtained in this way  (see Appendix~\ref{sec:A1}). The measure $\mu$ is called the exponent measure associated to the max-i.d. process $\eta$.

We now introduce the so-called $\alpha$- and $\beta$-mixing coefficients.  For more details on strong mixing conditions, the reader should refer to the recent survey by C.Bradley \cite{B05} or to the monographs \cite{D94,R00,B07a,B07b,B07c,DD07}.
For  $S\subset T$ a  closed subset, we denote by $\cF_S$ the $\sigma$-field generated by the random variables $\{\eta(s),\ s\in S\}$. Let $S_1,S_2\subset T$ be disjoint closed subsets. The $\alpha$-mixing coefficient (or strong mixing coefficient by Rosenblatt \cite{R56})  between the $\sigma$-fields $\cF_{S_1}$ and $\cF_{S_2}$  is defined by
\[
\alpha(S_1,S_2)=\sup\Big\{|\bbP(A\cap B)-\bbP(A)\bbP(B)|;\ A\in \cF_{S_1}, B\in \cF_{S_2} \Big\}.
\]
The $\beta$-mixing coefficient (or absolute regularity coefficient, see Volkonskii and Rozanov \cite{VR59}) between the $\sigma$-fields $\cF_{S_1}$ and $\cF_{S_2}$  is given by
\begin{eqnarray}
\beta(S_1,S_2)&=&\| \cP_{S_1\cup S_2}-\cP_{S_1}\otimes \cP_{S_2}\|_{var}\label{eq:bmix}\\
&=& \sup\Big\{ |\cP_{S_1\cup S_2}(C)-\cP_{S_1}\otimes \cP_{S_2}(C)|;\ C\in \cC_{S_1\cup S_2}\Big\}\nonumber
\end{eqnarray}
where $\|\cdot\|_{var}$ denotes the total variation of a signed measure and $\cP_S$ is the distribution of the restriction $\eta_{S}$ in the set $\bbC(S)$ of continuous functions on $S$ endowed with its Borel $\sigma$-field $\cC_S$. Since  $S_1$ and $S_2$ are disjoint closed subsets,  $\bbC(S_1\cup S_2)$ is canonically identified with $\bbC(S_1)\times\bbC(S_2)$.
It is well known that
\begin{eqnarray*}
\beta(S_1,S_2)&=& \frac{1}{2}\sup\Big\{ \sum_{i=1}^I\sum_{j=1}^J |\bbP(A_i\cap B_j)-\bbP(A_i)\bbP(B_j)|\Big\}
\end{eqnarray*}
where the supremum is taken over all partitions $\{A_1,\ldots,A_I\}$ and $\{B_1,\ldots,B_J\}$ of $\Omega$ with the $A_i$'s in $\cF_{S_1}$ and the $B_j$'s in $\cF_{S_2}$. It also holds 
\begin{equation}\label{eq:abmix}
\alpha(S_1,S_2)\leq \frac{1}{2}\beta(S_1,S_2). 
\end{equation}

Our main result is the following.
\begin{theo}\label{theo1}
Let $\eta$ be a max-i.d. process on $T$ with exponent measure $\mu$. Then, for all disjoint closed subsets $S_1,S_2\subset T$,
 \[
  \beta(S_1,S_2)\leq 2 \int_{\bbsC_0}\bbP[f\not<_{S_1}\eta,\ f\not<_{S_2}\eta]\,\mu(df).
 \]
\end{theo}
In the particular case when   $S_1$ and $S_2$ are finite or countable (which naturally arise for example if $T=\bbZ^d$), we can provide an upper bound for the mixing coefficient $\beta(S_1,S_2)$ involving only the $2$-dimensional marginal distributions of the process $\eta$. 

For $(s_1,s_2)\in T^2$, let $\mu_{s_1,s_2}$ be the exponent measure of the max-i.d. random vector $(\eta(s_1),\eta(s_2))$ defined  on $[0,+\infty)^2$ by
\[
\mu_{s_1,s_2}(A)=\mu\big[\{f\in\bbC_0(T);\ (f(s_1),f(s_2))\in A\}\big],\quad A\subset[0,+\infty)^2 \ \mbox{Borel\ set}.
\]

\begin{cor}\label{cor1}
If $S_1$ and $S_2$ are finite or countable disjoint closed subsets of $T$,
 \[
  \beta(S_1,S_2)\leq 2\sum_{s_1\in S_1}\sum_{s_2\in S_2} \int\bbP[\eta(s_1)\leq y_1,\ \eta(s_2)\leq y_2]\,\mu_{s_1,s_2}(dy_1dy_2).
 \]
\end{cor}

Next, we focus on simple max-stable random fields, where the phrase simple means that the marginals are standardized to the standard unit Fr\'echet distribution, 
\[
\bbP[\eta(t)\leq y]=\exp[-y^{-1}]1_{\{y>0\}},\quad y\in \bbR,\ t\in T.
\] 
In this framework, an insight into the dependence structure is given by the  extremal coefficients $\theta(S)$, $S\subset T$ compact, defined by the relation
\begin{equation}\label{eq:extcoeff}
\bbP[\sup_{s\in S} \eta(s)\leq y]=\exp[-\theta(S) y^{-1}],\quad y>0.
\end{equation}

\begin{theo}\label{theo2}
Let $\eta$ be a simple max-stable random field.\\
For all compact $S\subset T$, 
\begin{equation}\label{eq:defC}
C(S)=\bbE[\sup\{\eta(s)^{-1};\ s\in S\}]\ <\infty
\end{equation}
and furthermore:
\begin{itemize}
 \item For all disjoint compact  subsets $S_1,S_2\subset T$,
 \[
  \beta(S_1,S_2)\leq 2\Big[C(S_1)+C(S_2)\Big]  \Big[\theta(S_1)+\theta(S_2)-\theta(S_1\cup S_2)\Big].
 \]
\item For  $(S_{1,i})_{i\in I}$ and $(S_{2,j})_{j\in J}$  countable families of compact subsets of $T$ such that $S_1=\cup_{i\in I} S_{1,i}$ and $S_2=\cup_{j\in J} S_{2,j}$ are disjoint,
\[
  \beta(S_1,S_2)\leq 2\sum_{ i\in I}\sum_{ j\in J}\Big[C(S_{1,i})+C(S_{2,j})\Big]  \Big[\theta(S_{1,i})+\theta(S_{2,j})-\theta(S_{1,i}\cup S_{2,j})\Big].
 \]
\end{itemize}
\end{theo}
In the particular case when $S_1$ and $S_2$ are finite or countable, the mixing coefficient $\beta(S_1,S_2)$ can be bounded from above in terms of the extremal coefficient function
\[
\theta(s_1,s_2)=\theta(\{s_1,s_2\}),\quad s_1,s_2\in T. 
\]
 We recall the following basic properties:
it always holds $\theta(s_1,s_2)\in [1,2]$;
 $\theta(s_1,s_2)=2$ iff $\eta(s_1)$ and $\eta(s_2)$ are independent;
 $\theta(s_1,s_2)=1$ iff $\eta(s_1)=\eta(s_2)$.
Thus the extremal coefficient function gives some insight into the $2$-dimensional dependence structure of the max-stable field $\eta$, although it does not characterize it completely.

\begin{cor}\label{cor2}
Suppose $\eta$ is a simple max-stable random fields. If $S_1$ and $S_2$ are finite or countable disjoint closed subsets of $T$, then
 \[
  \beta(S_1,S_2)\leq 4\sum_{s_1\in S_1}\sum_{s_2\in S_2} [2-\theta(s_1,s_2)].
 \]
\end{cor}
We illustrate our results on two classes of stationary max-stable random fields on $\bbR^d$.
\begin{ex}\rm{ We consider the Brown-Resnick simple max-stable model (see Kabluchko et al. \cite{KSdH09}). Let $(W_i)_{i\geq 1}$ be independent copies of a  sample continuous  stationary increments Gaussian random field $W=(W(t))_{t\in\bbsR^d}$ with zero mean and variance $\sigma^2(t)$. Independently, let $(Z_i)_{i\geq 1}$ be the nonincreasing enumeration of the points of a Poisson point process  $(0,+\infty)$ with  intensity $z^{-2}dz$. The associated Brown-Resnick max-stable random field  is defined by
\[
\eta(t)=\bigvee_{i=1}^\infty Z_i \exp[W_i(t)-\sigma^2(t)/2],\quad t\in\bbR^d.
\]
It is known that $\eta$ is a stationary simple max-stable random fields whose law depends only on the  semi-definite negative function $V$, called the variogram of $W$, and defined by
\[
V(h)=\bbE[(W(t+h)-W(t))^2],\quad h\in\bbR^d.
\]
In this case, the extremal coefficient function is given by
\[
\theta(s_1,s_2)= 2\Psi(\sqrt{V(s_2-s_1)}/2),\quad s_1,s_2\in\bbR^d,
\]
where $\Psi$ denotes the cdf of the standard normal law.

}
\end{ex}
\begin{ex}\rm{Our second class of example is the moving maximum process by de Haan and Pereira \cite{dHP06}. Let $f:\bbR^d\to[0,+\infty)$ be a continuous density function such that
\[
\int_{\bbsR^d}f(x)\,dx=1\quad \mbox{and}\quad \int_{\bbsR^d}\sup_{|h|\leq 1} f(x+h)\,dx<\infty.
\]
Let $\sum_{i=1}^\infty\delta_{(Z_i,U_i)}$  be a Poisson random measure on $(0,+\infty)\times\bbR^d$ with intensity $z^{-2}dzdu$.
Then the random field
\[
\eta(t)=\bigvee_{i=1}^\infty Z_i f(t-U_i),\quad t\in\bbR^d,
\]
is a stationary sample continuous simple max-stable random field. The corresponding extremal coefficient function is given by
\[
\theta(s_1,s_2)=\int_{\bbsR^d}\max(f(s_1-x),f(s_2-x))\,dx,\quad s_1,s_2\in\bbR^d.
\]
More generally, for any compact $S\subset \bbR^d$,
\[
\theta(S)=\int_{\bbsR^d}\max_{s\in S}f(s-x)\,dx.
\]
}
\end{ex}

As noted in the introduction, our main  motivation for considering the strong mixing properties of max-i.d. random fields is to obtain central limit Theorems (CLTs) for stationary max-i.d. random fields. In this direction, we focus on stationary random fields on $T=\bbZ^d$ and our analysis relies on  Bolthausen's CLT  \cite{B82} (see Appendix~\ref{sec:TCL}). 

We note $|h|=\max_{1\leq i\leq d} |h_i|$ the norm of $h\in\bbZ^d$ and   $|S|$ the number of elements of a subset $S\subset \bbZ^d$. We note $\partial S$ the set of elements $h\in S$ such that there is $h'\notin S$ with $d(h,h')=1$. The distance between two subsets $S_1,S_2\subset \bbZ^d$ is given by 
$$d(S_1,S_2)=\min\{|s_2-s_1|;\ s_1\in S_1,\ s_2\in S_2\}.$$ 

A random field $X=(X(t))_{t\in\bbsZ^d}$ is said to be stationary if the law of $(X(t+s))_{t\in\bbsZ^d}$ does not depend on $s\in\bbZ^d$. We say that a square integrable stationary random field $X$ satisfies the CLT if the following two conditions are satisfied:
\begin{itemize}
\item[i)] the series $\sigma^2=\sum_{t\in\bbsZ^d} \mathrm{Cov}[X(0),X(t)]$ converges absolutely;
\item[ii)] for all  sequence  $\Lambda_n$ of  finite subsets of $\bbZ^d$, which increases to $\bbZ^d$ and such that $\lim_{n\to\infty}|\partial\Lambda_n|/|\Lambda_n|=0$, the sequence
$|\Lambda_n|^{-1/2}\sum_{t\in\Lambda_n} (X(t)-\bbE[X(t)])$
converge in law to the normal distribution with mean $0$ and variance $\sigma^2$ as $n\to\infty$.
\end{itemize}
Please note we do not require the limit variance $\sigma^2$ to be positive; the case $\sigma^2=0$ corresponds to a degenerated CLT where the limit distribution is the Dirac mass at zero. Bolthausen's CLT for stationary mixing random fields together with our estimates of  mixing coefficients of max-i.d. random fields yield the following Theorem.
\begin{theo}\label{theo3}
Suppose $\eta$ is a stationary max-i.d. random fields on $\bbZ^d$ with exponent measure $\mu$ and let
\[
\gamma(h)=\int\bbP[\eta(0)\leq y_1,\ \eta(h)\leq y_2]\,\mu_{0,h}(dy_1dy_2),\quad h\in\bbZ^d.
\]
Let  $g:\bbR^p\to \bbR$ be a measurable function and $t_1,\ldots,t_p\in \bbZ^d$ such that 
\[
\bbE[g(\eta(t_1),\ldots,\eta(t_p))^{2+\delta}]<\infty \ \mbox{ for some }\ \delta>0,
\]
and assume that 
\begin{equation}\label{eq:condgamma}
\sum_{|h|\geq m}\gamma(h)=o(m^{d-1}) \quad \mbox{and}\quad  \sum_{m=1}^\infty m^{d-1}\sup_{|h|\geq m}\gamma(h)^{\delta/(2+\delta)}<\infty.
\end{equation}
Then the stationary random field $X$ defined by
\[
X(t)=g(\eta(t_1+t),\ldots,\eta(t_p+t)),\quad t\in\bbZ^d
\]
 satisfies the CLT.
\end{theo}
Condition \eqref{eq:condgamma} requires that $\gamma$ goes fast enough to $0$ at infinity. It is met for example if 
\begin{equation}\label{eq:condgamma2}
\gamma(h)\leq C\cdot|h|^{-b}\quad \mbox{for\ some\ }\ \ b>d\max(2,(2+\delta)/\delta)\ \ \mbox{and}\ \ C>0.
\end{equation}
If $\eta$ is simple max-stable, we have $\gamma(h)\leq 2(2-\theta(0,h))$, with $\theta$ the extremal coefficient function.

As an application, we consider the estimation of the extremal coefficient for a stationary simple max-stable random field on $\bbZ^d$. For $h\in\bbZ^d$, we note $\theta(h)=\theta(0,h)$. Equation \eqref{eq:extcoeff} implies
\[
\theta(h)=-y\log p(h,y)\quad \mbox{with}\quad p(h,y)=\bbP(\eta(0)\leq y,\ \eta(h)\leq y),\quad y>0,
\]
suggesting the simple estimator
\begin{equation}\label{eq:hattheta}
\hat\theta_n^{(1)}(h)=-y\log \hat p_n(h,y)
\quad \mbox{with}\quad  \hat p_n(h,y)= |\Lambda_n|^{-1}\sum_{t\in\Lambda_n}1_{\{\eta(t)\leq y,\ \eta(t+h)\leq y\}}
\end{equation}
where $\Lambda_n$ is a sequence of finite subsets increasing to $\bbZ^d$  such that $|\partial\Lambda_n|/|\Lambda_n|\to 0$ as $n\to\infty$. The fact that the naive estimator $\hat\theta_n^{(1)}(h)$ depends of the threshold level $y>0$ is not satisfactory. Alternatively, one may consider the following procedures. Smith \cite{S90} noticed that $\min(\eta(0)^{-1},\eta(h)^{-1})$ as an exponential distribution with mean $\theta(h)^{-1}$ and proposed the estimator
\[
\hat\theta^{(2)}_n(h)=\frac{|\Lambda_n|}{\sum_{t\in\Lambda_n} \min(\eta(t)^{-1},\eta(t+h)^{-1})}.
\]
Cooley et al. \cite{CNP06} introduced the $F$-madogram defined by
\[
\nu_F(h)=\bbE[|F(\eta(0))-F(\eta(h))|]\quad \mbox{with}\quad F(y)=\exp(-1/y)1_{\{y>0\}}
\]
and show that it satisfies
\[
\nu_F(h)=\frac{1}{2}\frac{\theta(h)-1}{\theta(h)+1}\quad \mbox{or\ equivalently}\quad \theta(h)=\frac{1+2\nu_F(h)}{1-2\nu_F(h)}.
\]
This suggests the estimator
\[
\hat\theta_n^{(3)}(h)=\frac{|\Lambda_n|+2\sum_{t\in\Lambda_n} |F(\eta(t))-F(\eta(t+h))|}{|\Lambda_n|-2\sum_{t\in\Lambda_n} |F(\eta(t))-F(\eta(t+h))|}.
\]
The following Proposition states the asymptotic normality of these estimators. 
\begin{prop}\label{prop:asynorm}
Suppose that $\eta$ is a stationary random field on $\bbZ^d$ with extremal coefficient function satisfying 
\begin{equation}\label{eq:asynorm1}
2-\theta(h)\leq C\cdot|h|^{-b}\quad \mbox{for\ some\ }\ \ b>2d\ \ \mbox{and}\ \ C>0.
\end{equation}
Then, the estimators $\hat\theta_n^{(i)}(h)$, $i=1,2,3$ are asymptotically normal: 
\begin{eqnarray*}
|\Lambda_n|^{1/2}\big(\hat\theta_n^{(i)}(h)-\theta(h)\big)&\Longrightarrow& \cN(0,\sigma_i^2) \quad \mbox{as}\ n\to\infty
\end{eqnarray*}
with limit variances 
\begin{eqnarray*}
\sigma^{2}_1&=&y^2\sum_{t\in\bbsZ^d}\big(\exp[(2\theta(h)-\theta(\{0,h,t,t+h\}))y^{-1}]-1\big),\\
\sigma^{2}_2&=&\theta(h)^4\sum_{t\in\bbsZ^d}\mathrm{Cov}\big[\min(\eta(0)^{-1},\eta(h)^{-1})\,,\, \min(\eta(t)^{-1},\eta(t+h)^{-1})\big],\\
\sigma^{2}_3&=&(\theta(h)+1)^4\sum_{t\in\bbsZ^d}\mathrm{Cov}\big[|F(\eta(0))-F(\eta(h))|\,,\, |F(\eta(t))-F(\eta(t+h))|\big].\\
\end{eqnarray*}
\end{prop}
Interestingly, the function $y\mapsto \sigma^2_1$ is strictly convex, has limit $+\infty$ as $y\to 0^+$ or $+\infty$ and hence it admits a unique minimizer  $y^\star$ corresponding to an asymptotically optimal truncation level for the estimator $\hat\theta_n^{(1)}$. Unfortunately, the limit variances $\sigma_2^2$ and $\sigma_3^2$ are  not so explicit so that a comparison between the three is difficult.

\section{Proofs}
\subsection{Strong mixing properties of extremal point processes}
In the sequel, we shall write shortly $\bbC_0=\bbC_0(T)$. We introduce here the notion of $S$-extremal points that will play a key role in this work. We use the following notations: if $f_1,f_2$ are two functions
defined (at least) on $S$, we note
\begin{align*}
  f_1=_Sf_2\quad&\text{if and only if}\quad \forall s\in S,\ f_1(s)=f_2(s),\\
  f_1<_Sf_2\quad&\text{if and only if}\quad \forall s\in S,\   f_1(s)<f_2(s),\\
  f_1\not<_Sf_2\quad&\text{if and only if}\quad \exists s\in S,\ f_1(s)\geq f_2(s).
\end{align*}
A point $\phi\in\Phi$ is said to be $S$-subextremal if $\phi<_S \eta$, it is said $S$-extremal otherwise, i.e. if there exists $s\in S$ such that $\phi(s)=\eta(s)$. In words, a $S$-subextremal point has no contribution to the maximum $\eta$ on $S$.
\begin{df}\label{def:ep}
Define the $S$-extremal random point process $\Phi_S^+$ and the $S$-subextremal random point process $\Phi_S^-$ by
\[
\Phi_S^+=\{\phi\in\Phi;\ \phi\not<_S\eta\}\quad \mbox{and}\quad \Phi_S^-=\{\phi\in\Phi;\ \phi<_S\eta\}.
\]
\end{df}
\noindent
The fact that $\Phi_S^+$ and $\Phi_S^-$ are well defined point processes, i.e. that they satisfy some measurability properties, is proved in \cite{DEM11b} Appendix A3.  Clearly, the restriction  $\eta_{S}$  depends on $\Phi_S^+$ only:
\[
\eta(s)=\max\{\phi(s);\ \phi\in \Phi_S^+\},\quad s\in S.
\]
This implies that the  strong mixing coefficient $\beta(S_1,S_2)$ defined by Equation \eqref{eq:bmix} can be upper bounded by a similar $\beta$-mixing coefficient defined on the level of the extremal point process $\Phi_{S_1}^+$, $\Phi_{S_2}^+$. For $i=1,2$, let $P_{\Phi_{S_i}^+}$ the distribution of $\Phi_{S_i}^+$ on the space of locally finite point measures on $\bbC_0$ and let  $P_{(\Phi_{S_1}^+,\Phi_{S_2}^+)}$ be the joint distribution of $(\Phi_{S_1}^+,\Phi_{S_2}^+)$. We define
\begin{equation}\label{eq:bmix2}
\beta(\Phi_{S_1}^+,\Phi_{S_2}^+)=\|P_{(\Phi_{S_1}^+,\Phi_{S_2}^+)}-P_{\Phi_{S_1}^+}\otimes P_{\Phi_{S_2}^+}\|_{var}.
\end{equation}
It holds
\begin{equation}\label{eq:comp}
\beta(S_1,S_2)\leq \beta(\Phi_{S_1}^+,\Phi_{S_2}^+). 
\end{equation}
The following Theorem provides a simple estimate for the $\beta$-mixing coefficient on the point process level. It implies Theorem \ref{theo1}  straightforwardly and has a clearer interpretation.
\begin{theo}\label{theo4}
\begin{itemize}
\item The following upper bound holds true:
\begin{equation}\label{eq:bmix4}
\beta(\Phi_{S_1}^+,\Phi_{S_2}^+)\leq 2\,\bbP[\Phi_{S_1}^+\cap \Phi_{S_2}^+= \emptyset]
\end{equation}
with
\begin{equation}\label{eq:bmix4bis}
\bbP[\Phi_{S_1}^+\cap \Phi_{S_2}^+= \emptyset]\leq \int_{\bbsC_0}\bbP(f\not<_{S_1}\eta,\ f\not<_{S_2}\eta)\,\mu(df).
\end{equation}
\item If the point process $\Phi$ is simple (in particular in the max-stable case), the following lower bound holds true:
\begin{equation}\label{eq:bmix4ter}
\beta(\Phi_{S_1}^+,\Phi_{S_2}^+)\geq \bbP[\Phi_{S_1}^+\cap \Phi_{S_2}^+= \emptyset]
\end{equation}
\end{itemize}
\end{theo}

Clearly, equations \eqref{eq:comp}, \eqref{eq:bmix4} and \eqref{eq:bmix4bis} together imply Theorem \ref{theo1}.

\begin{rem}{\rm The upper and lower bound in Theorem \ref{theo4} are of the same order, and hence relatively sharp. It is not clear however how to lower bound  $\beta(S_1,S_2)$ and  how sharp the upper bound in Theorem \ref{theo1} is.}
\end{rem}

\subsection{Proof of Theorem \ref{theo4}}

The upper bound for the mixing coefficient $\beta(\Phi_{S_1}^+,\Phi_{S_2}^+)$ defined by Equation \eqref{eq:bmix2} relies on a standard coupling argument. There are indeed deep relationships between $\beta$-mixing  and optimal couplings (see e.g. \cite{R00} chapter 5) and we will use the following well known result:
\begin{lem}\label{lem:coupling}
On a probability space $(\Omega,\cF,\bbP)$, suppose that the random variables \\ $(\Phi_{S_1}^{+i},\Phi_{S_2}^{+i})$, $i=1,2$ are  such that:
\begin{itemize}
\item[i)] the distribution of $(\Phi_{S_1}^{+1},\Phi_{S_2}^{+1})$ is $P_{(\Phi_{S_1}^+,\Phi_{S_2}^+)}$;
\item[ii)] the distribution of $(\Phi_{S_1}^{+2},\Phi_{S_2}^{+2})$ is $P_{\Phi_{S_1}^+}\otimes P_{\Phi_{S_2}^+}$.
\end{itemize}
Then, 
\[
\beta( \Phi_{S_1}^+,\Phi_{S_2}^+)\leq \bbP\Big[(\Phi_{S_1}^{+1},\Phi_{S_2}^{+1})\neq (\Phi_{S_1}^{+2},\Phi_{S_2}^{+2})\Big]
\]
\end{lem}
We say that the random variables $(\Phi_{S_1}^{+i},\Phi_{S_2}^{+i})$, $i=1,2$ satisfying i) and ii) realize a coupling between the distributions $P_{(\Phi_{S_1}^+,\Phi_{S_2}^+)}$ and $P_{\Phi_{S_1}^+}\otimes P_{\Phi_{S_2}^+}$. 

In order to construct a suitable coupling, we need the following Lemma describing  the dependence between $\Phi_S^+$ and $\Phi_S^-$. 
\begin{lem}\label{lem:lc}
Let $S\subset T$ be a closed set. The conditional distribution of $\Phi_S^-$ with respect to $\Phi_S^+$ is equal to the distribution of a Poisson point process  with intensity $1_{\{f<_S \eta\}}\mu(df)$.
\end{lem}
Note that in the particular case when $T$ is compact and $\Phi_S^+$ is finite almost surely, Lemma 
\ref{lem:lc} follows from  \cite{DEM11b} Theorem 2.1 and Corollary 2.1. For $T$ non compact, the proof need to be modified in a non straightforward way and is postponed to the Appendix \ref{sec:lemlc}.

We now construct the coupling providing the upper bound \eqref{eq:bmix4}.
\begin{prop}\label{prop:coupling}
Let $(\widetilde\Phi,\widetilde\eta)$ be an independent copy of $(\Phi,\eta)$
and define 
\begin{equation}\label{eq:defhatphi}
\widehat\Phi=\Phi_{S_1}^+ \cup \{\widetilde \phi\in \widetilde\Phi; \ \widetilde\phi<_{S_1} \eta\}.
\end{equation}
The following properties hold true:
\begin{itemize}
\item  $\widehat\Phi$ has the same distribution as $\Phi$ and $\widetilde\Phi$ and satisfies 
\begin{equation}\label{eq:dec}
\widehat\Phi_{S_1}^+=\Phi_{S_1}^+\quad ,\quad \widehat\Phi_{S_1}^-=\{\widetilde \phi\in \widetilde\Phi; \ \widetilde\phi<_{S_1} \eta\};
\end{equation}
\item $(\widehat\Phi_{S_1}^+,\widehat\Phi_{S_2}^+)$ and $(\Phi_{S_1}^+,\widetilde \Phi_{S_2}^+)$ is a coupling between $P_{(\Phi_{S_1}^+,\Phi_{S_2}^+)}$ and $P_{\Phi_{S_1}^+}\otimes P_{\Phi_{S_2}^+}$ such that
\begin{equation}\label{eq:coupling}
\bbP\Big[(\widehat\Phi_{S_1}^+,\widehat\Phi_{S_2}^+)\neq (\Phi_{S_1}^+,\widetilde\Phi_{S_2}^+) \Big]
\leq 2\,\bbP[\Phi_{S_1}^+\cap\Phi_{S_2}^+=\emptyset] . \end{equation}
\end{itemize}
\end{prop}

\begin{preuve}{of Proposition \ref{prop:coupling}}
\begin{itemize}
\item Equation \eqref{eq:dec} follows from the construction of $\widehat\Phi$: consider
\[
\widehat\eta(t)=\bigvee_{\phi\in\widehat\Phi}\phi(t),\quad t\in T;
\]
the maximum $\eta$ is achieved on $S_1$ by the $S_1$-extremal points $\Phi_{S_1}^+$, and the definition \eqref{eq:defhatphi} ensures that $\widehat\eta$ and $\eta$ are equal on $S_1$ so that equation \eqref{eq:dec} holds.

Furthermore, conditionally on $\Phi_{S_1}^+$, the distribution of $\{\widetilde \phi\in \widetilde\Phi; \ \widetilde\phi<_{S_1} \eta\}$ is  equal to the distribution of a Poisson point process with intensity $1_{\{f<_{S_1}\eta\}}\mu(df)$. 
According to Lemma \ref{lem:lc}, this is the conditional distribution of $\Phi_{S_1}^-$ given $\Phi_{S_1}^+$, whence $(\widehat\Phi_{S_1}^+,\widehat\Phi_{S_1}^-)$  has the same distribution as $(\Phi_{S_1}^+,\Phi_{S_1}^-)$. We deduce that $\Phi=\Phi_{S_1}^+\cup\Phi_{S_1}^-$ and $\widehat\Phi=\widehat\Phi_{S_1}^+\cup\widehat\Phi_{S_1}^-$ have the same distribution.

\item The coupling property is easily proved: since $\Phi$ and $\widehat\Phi$ have the same distribution, the law of $(\widehat\Phi_{S_1}^+,\widehat\Phi_{S_2}^+)$ is equal to $P_{(\Phi_{S_1}^+,\Phi_{S_2}^+)}$;  since $\Phi$ and $\widetilde \Phi$ are independent, $(\Phi_{S_1}^+,\widetilde \Phi_{S_2}^+)$ has law $P_{\Phi_{S_1}^+}\otimes P_{\Phi_{S_2}^+}$.\\ 
We are left to prove Equation \eqref{eq:coupling}. Since $\Phi_{S_1}^+=\widehat\Phi_{S_1}^+$, we need to upper bound the probability
$\bbP\big[\widehat\Phi_{S_2}^+\neq \widetilde\Phi_{S_2}^+\big]$. By construction, $\widehat\Phi$ is obtained from $\widetilde\Phi$ by removing the points $\widetilde\phi\in\widetilde\Phi$ such that $\widetilde\phi\not <_{S_1} \eta$ and  adding the points $\phi\in \Phi_{S_1}^+$. Hence, it holds
\[
 \{\widehat \Phi_{S_2}^+ \neq  \widetilde  \Phi_{S_2}^+\} \ \subset \ \{\exists \phi\in \Phi_{S_1}^+,\ \phi\not<_{S_2} \widehat \eta  \}\  \cup\  \{  \exists \widetilde \phi\in \widetilde \Phi_{S_2}^+,\ \widetilde\phi \not<_{S_1}  \eta \}
\]
Noting the equality of events  
\[
\{\exists \phi\in \Phi_{S_1}^+,\ \phi\not<_{S_2} \widehat \eta  \}=\{\exists \phi\in \widehat\Phi_{S_1}^+\cap \widehat\Phi_{S_2}^+\}=\{\widehat\Phi_{S_1}^+\cap\widehat\Phi_{S_2}^+\neq \emptyset\},
\]
we obtain
\[
\bbP[\widehat \Phi_{S_2}^+ \neq  \widetilde  \Phi_{S_2}^+]\leq\bbP[\widehat\Phi_{S_1}^+\cap\widehat\Phi_{S_2}^+\neq \emptyset]+\bbP[\exists \widetilde \phi\in \widetilde \Phi_{S_2}^+,\ \widetilde\phi \not<_{S_1}  \eta ].
\]
Since $\widehat\Phi$ and $\Phi$ have the same law, we have 
\[
\bbP[\widehat\Phi_{S_1}^+\cap\widehat\Phi_{S_2}^+\neq \emptyset]=\bbP[\Phi_{S_1}^+\cap\Phi_{S_2}^+\neq \emptyset].
\]
Hence, equation \eqref{eq:coupling} follows from the upper bound
\[
\bbP[\exists \widetilde \phi\in \widetilde \Phi_{S_2}^+,\ \widetilde\phi \not<_{S_1}  \eta ]\leq \bbP[\Phi_{S_1}^+\cap\Phi_{S_2}^+\neq \emptyset]
\]
that we prove now.
Using  symmetry and exchanging the roles of $\Phi$ and $\widetilde\Phi$ on the one hand and the roles of $S_1$ and $S_2$ on the other hand, it is equivalent to prove that
\[
\bbP[\exists  \phi\in  \Phi_{S_1}^+,\ \phi \not<_{S_2}  \widetilde\eta ]\leq \bbP[\Phi_{S_1}^+\cap\Phi_{S_2}^+\neq \emptyset].
\]
We conclude the proof by noticing that the inclusion of events
\[
\{\exists  \phi\in  \Phi_{S_1}^+,\ \phi \not<_{S_2}  \widetilde\eta \}\subset \{\exists \phi\in \Phi_{S_1}^+,\ \phi\not<_{S_2} \widehat \eta  \}
\]
entails
\[
\bbP[\exists  \phi\in  \Phi_{S_1}^+,\ \phi \not<_{S_2}  \widetilde\eta ]\leq \bbP[\exists  \phi\in  \Phi_{S_1}^+,\ \phi \not<_{S_2}  \widehat\eta ]=\bbP[\Phi_{S_1}^+\cap\Phi_{S_2}^+\neq \emptyset].
\]
\end{itemize}
\end{preuve}

We now complete the proof of Theorem \ref{theo4} by proving Equations  \eqref{eq:bmix4bis} and \eqref{eq:bmix4ter}.

\begin{preuve}{of Equation \eqref{eq:bmix4bis}}
We observe that
\[
\{\Phi_{S_1}^+\cap \Phi_{S_2}^+\neq\emptyset\}=\{\exists\phi\in\Phi,\ \phi\not<_{S_1}\eta,\ \phi\not<_{S_1}\eta\}
\]
which entails
\[
\bbP[\Phi_{S_1}^+\cap \Phi_{S_2}^+\neq\emptyset]
\leq \bbE\Big[\sum_{\phi\in\Phi} 1_{\{\phi\not<_{S_1}\eta,\ \phi\not<_{S_2}\eta \}}\Big].
\]
Noting that $\phi\not<_{S_i}\eta$ if and only if $\phi\not<_{S_i} \max(\Phi- \{\phi\})$, we apply Slyvniak's formula (see Appendix ) and compute
\begin{eqnarray*}
\bbE\Big[\sum_{\phi\in\Phi} 1_{\{\phi\not<_{S_1}\eta,\ \phi\not<_{S_2}\eta \}}\Big]
&=&\bbE\Big[\sum_{\phi\in\Phi} 1_{\{\phi\not<_{S_1}\max(\Phi- \{\phi\}),\ \phi\not<_{S_2}\max(\Phi- \{\phi\}) \}}\Big]\\
&=& \int_{\bbsC_0}\bbE[1_{\{f \not<_{S_1}\max(\Phi),\ f \not<_{S_2}\max(\Phi)\}}]\mu(df)\\
&=& \int_{\bbsC_0}\bbP[f \not<_{S_1}\eta,\ f \not<_{S_2}\eta]\mu(df).
\end{eqnarray*}
\end{preuve}

\begin{preuve}{of Equation \eqref{eq:bmix4ter}}
Clearly, for all measurable subset $C\subset M_p(\bbC_0)\times M_p(\bbC_0)$,
\[
\beta(\Phi_{S_1}^+,\Phi_{S_2}^+)\geq |\bbP[(\Phi_{S_1}^+,\Phi_{S_2}^+)\in C]-\bbP[(\Phi_{S_1}^+,\widetilde\Phi_{S_2}^+)\in C]|
\] 
where $\widetilde\Phi$ is an independent copy of $\Phi$. We obtain the lower bound \eqref{eq:bmix4ter} by choosing the subset 
\[
C=\{(M_1,M_2)\in M_p(\bbC_0)\times M_p(\bbC_0);\ M_1\cap M_2\neq \emptyset\}.
\]
This yields indeed 
\[
\beta(\Phi_{S_1}^+,\Phi_{S_2}^+)\geq |\bbP[\Phi_{S_1}^+\cap \Phi_{S_2}^+\neq \emptyset]-\bbP[\Phi_{S_1}^+\cap \widetilde\Phi_{S_2}^+\neq \emptyset]|.
\] 
and in the case when $\Phi$ is a simple point process, i.e. when the intensity measure $\mu$ has no atom, we have
\[
\bbP[\Phi_{S_1}^+\cap\widetilde \Phi_{S_2}^+\neq\emptyset]= \bbP[\Phi\cap\widetilde \Phi\neq\emptyset]=0.
\]
\end{preuve}

\subsection{Proof of Corollaries \ref{cor1} and \ref{cor2} and Theorem \ref{theo2}}

\begin{preuve}{of Corollary \ref{cor1}} We have for all $f\in\bbC_0$,
\begin{eqnarray*}
  \{f\not<_{S_1}\eta,\ f\not<_{S_2}\eta\}
&=& \{\exists(s_1,s_2)\in S_1\times S_2,\ f(s_1)\geq \eta(s_1),\ f(s_2)\geq \eta(s_2)\}\\
&=&\cup_{s_1\in S_1}\cup_{s_2\in S_2}\{\eta(s_1)\leq f(s_1),\ \eta(s_2)\leq f(s_2)\}
\end{eqnarray*}
whence, for  $S_1$ and $S_2$  finite  or countable, 
\begin{eqnarray*}
&&\bbP[f\not<_{S_1}\eta,\ f\not<_{S_2}\eta]\leq  \sum_{s_1\in S_1}\sum_{s_2\in S_2}\bbP[\eta(s_1)\leq f(s_1),\ \eta(s_2)\leq f(s_2)].
\end{eqnarray*}
As a consequence, the integral in Theorem \ref{theo1} satisfies
\begin{eqnarray*}
&&\int_{\bbsC_0}\bbP[f\not<_{S_1}\eta,\ f\not<_{S_2}\eta]\, \mu(df)\\
&\leq & \sum_{s_1\in S_1}\sum_{s_2\in S_2}\int_{\bbsC_0}\bbP[\eta(s_1)\leq f(s_1),\ \eta(s_2)\leq f(s_2)]\,\mu(df)\\
&=&\sum_{s_1\in S_1}\sum_{s_2\in S_2}\int_{[0,+\infty)^2}\bbP[ \eta(s_1)\leq y_1,\  \eta(s_2)\leq y_2]\,\mu_{s_1,s_2}(dy_1dy_2).
\end{eqnarray*}
In the last line, we have used the fact that $\mu_{s_1,s_2}$ is the image of the measure $\mu$ under the mapping $f\mapsto (f(s_1),f(s_2))$.
\end{preuve}

\begin{preuve}{of Theorem \ref{theo2}}
We recall that for a simple max-stable random field, the exponent measure $\mu$ is homogeneous of order $-1$, i.e. $\mu(cA)=c^{-1}\mu(A)$ for all $A\subset \bbC_0$ Borel set and $c>0$. Also the assumption of standard unit Fréchet marginals implies
\[
\bar\mu_t(y)=y^{-1},\quad t\in T,\ y>0.
\]
These conditions implies (See Gin\'e and al. \cite{GHV90} Propostition 3.2  or de Haan and Fereira \cite{dHF06} Theorem 9.4.1 and Corollary 9.4.2) that $\mu$ can be written as
\[
\mu(A)=\int_0^\infty\int_{\bbsC_0} 1_{\{rf\in A\}}\,r^{-2}dr\sigma(df)
\]
where $\sigma$ is a probability measure on $\bbC_0$ such that 
\[
\int_{\bbsC_0} f(t)\,\sigma(df)=1\quad \mbox{for\ all\ }t\in T,
\]
and
\[
\int_{\bbsC_0} \sup_{s\in S}f(s)\,\sigma(df)<\infty\quad \mbox{for\ all\ compact}S\subset T.
\]
Using this, note that for all compact $S\subset T$ and $y>0$,
\begin{eqnarray*}
\bbP[\sup_{s\in S}f(s)\leq y]&=&\exp[-\mu(\{f\in\bbC_0;\ \sup_{s\in S}f(s)>y\})]\\
&=&\exp\Big[-\int_{\bbsC_0} 1_{\{\sup_{s\in S}rf(s)>y\}}\, r^{-2}dr\sigma(df)\Big]\\
&=&\exp\Big[-y^{-1}\int_{\bbsC_0}\sup_{s\in S}f(s)\,\sigma(df)\Big].
\end{eqnarray*}
It follows that the extremal coefficient $\theta(S)$ defined by \eqref{eq:extcoeff} is equal to 
\begin{equation}\label{eq:defec}
\theta(S)=\int_{\bbsC_0}\sup_{s\in S}f(s)\,\sigma(df).
\end{equation}
We now consider equation \eqref{eq:defC} and observe that for $S\subset T$ compact,
\[
 C(S)=\bbE\Big[(\inf_{s\in S}\eta(s))^{-1} \Big]
\]
so that we need to provide a lower bound for $\inf_{s\in S}\eta(s)$. To this aim, we remark that
\[
\inf_{s\in S} \eta(s)= \inf_{s\in S} \max_{\phi\in\Phi} \phi(s) \geq \max_{\phi\in\Phi} \inf_{s\in S} \phi(s).
\]
The right hand side is a random variable with unit Fréchet distribution since
\begin{eqnarray*}
\bbP[\max_{\phi\in\Phi} \inf_{s\in S} \phi(s)\leq y]
&=&  \bbP[\forall\phi\in\Phi,\  \inf_{s\in S} \phi(s)\leq y]\\
&=& \exp(-\mu(\{f\in\bbC_0;\ \inf_{s\in S} f(s)> y\}))
\end{eqnarray*}
and
\[
 \mu(\{f\in\bbC_0;\ \inf_{s\in S} f(s)> y\})=y^{-1}\int_{\bbsC_0}\inf_{s\in S}f(s)\, \sigma(df).
\]
Hence, if $\int_{\bbsC_0}\inf_{s\in S}f(s)\, \sigma(df)> 0$, we obtain
\[
 C(S)=\bbE\Big[(\inf_{s\in S}\eta(s))^{-1} \Big]\leq \bbE\Big[ (\max_{\phi\in\Phi} \inf_{s\in S} \phi(s))^{-1}\Big]=\Big(\int_{\bbsC_0}\inf_{s\in S}f(s)\, \sigma(df)\Big)^{-1}<\infty.
\]
For arbitrary compact $S\subset T$, we may however have $\int_{\bbsC_0}\inf_{s\in S}f(s)\, \sigma(df)= 0$. But, if $S=B(s_0,\varepsilon)$ is a closed ball with center $s_0$ and radius $\varepsilon$, the monotone convergence Theorem implies 
\[
\int_{\bbsC_0}\inf_{s\in B(s_0,\varepsilon)}f(s)\, \sigma(df)\to \int_{\bbsC_0}f(s_0)\, \sigma(df)=1\quad \mbox{as}\ \varepsilon\to 0,
\]
so that $\int_{\bbsC_0}\inf_{s\in B(s_0,\varepsilon_0)}f(s)\, \sigma(df)>0$ and $C(B(s_0,\varepsilon_0))<\infty$ for $\varepsilon_0$ small enough. The result for general $S$ follows by a compacity argument: there exist $s_1,\ldots,s_k$ and $\varepsilon_1,\ldots,\varepsilon_k$ such that $S\subset \cup_{i=1}^k B(s_i,\varepsilon_i)$. Hence,
\[
 \sup_{s\in S} \eta(s)^{-1} \leq  \max_{1\leq i\leq k}  \sup_{s\in B(s_i,\varepsilon_i)} \eta(s)^{-1} \leq \sum_{i=1}^k  \sup_{s\in B(s_i,\varepsilon_i)} \eta(s)^{-1}
\]
and 
\[
 C(S)=\bbE\Big [\sup_{s\in S} \eta(s)^{-1}\Big] \leq \sum_{i=1}^k \bbE\Big[\sup_{s\in B(s_i,\varepsilon_i)} \eta(s)^{-1}\Big]=\sum_{i=1}^k C(B(s_i,\varepsilon_i))<\infty.
\]
This proves equation \eqref{eq:defC}. 
\begin{itemize}
 \item The upper bound for $\beta(S_1,S_2)$ given by Theorem \ref{theo1} can be expressed as
\begin{eqnarray}
\beta(S_1,S_2)
&\leq&2\int_{\bbsC_0}\bbP[f\not<_{S_1}\eta,\ f\not<_{S_2}\eta]\, \mu(df)\nonumber\\
&=&2\int_0^\infty\int_{\bbsC_0}\bbP[rf\not<_{S_1}\eta,\ rf\not<_{S_2}\eta]\, r^{-2}dr\sigma(df)\nonumber\\
&=& 2\int_0^\infty\int_{\bbsC_0} \bbE\Big[ 1_{\{r\geq \inf_{s_1\in S_1}\eta(s_1)/f(s_1),\ r\geq \inf_{s_2\in S_2}\eta(s_2)/f(s_2)\}}\Big]\, r^{-2}dr\sigma(df)\nonumber\\
&=& 2\int_{\bbsC_0} \bbE\Big[ \max\big(\inf_{s_1\in S_1}\frac{\eta(s_1)}{f(s_1)},\ \inf_{s_2\in S_2}\frac{\eta(s_2)}{f(s_2)}\big)^{-1}\Big]\,\sigma(df)\label{eq:ub1}.
\end{eqnarray}
We then introduce the upper bound
\begin{eqnarray*}
&&\max\Big(\inf_{s_1\in S_1}\frac{\eta(s_1)}{f(s_1)},\ \inf_{s_2\in S_2}\frac{\eta(s_2)}{f(s_2)}\Big)^{-1}\\
&\leq& \max(\sup_{s_1\in S_1}\eta(s_1)^{-1},\sup_{s_2\in S_2}\eta(s_2)^{-1})\ \min(\sup_{s_1\in S_1}f(s_1),\sup_{s_2\in S_2}f(s_2))
\end{eqnarray*}
whence we deduce
\begin{eqnarray*}
&&\beta(S_1,S_2)\\
&=& 2\bbE\Big[ \max\big(\sup_{s_1\in S_1}\eta(s_1)^{-1},\sup_{s_2\in S_2}\eta(s_2)^{-1}\big)\Big]\ \int_{\bbsC_0}\min(\sup_{s_1\in S_1}f(s_1),\sup_{s_2\in S_2}f(s_2))\,\sigma(df)\\
&\leq &2\bbE\Big[\sup_{s_1\in S_1}\eta(s_1)^{-1}+\sup_{s_2\in S_2}\eta(s_2)^{-1}\Big]\ \int_{\bbsC_0}\min(\sup_{s_1\in S_1}f(s_1),\sup_{s_2\in S_2}f(s_2))\,\sigma(df)\\
&=& 2\big[C(S_1)+C(S_2)\big]\ \big[\theta(S_1)+\theta(S_2)-\theta(S_1\cup S_2)\big].
\end{eqnarray*}
In the last equality, we use Equation \eqref{eq:defC} defining $C(S)$ together with the following simple observation: in view of Equation \eqref{eq:defec}, the equality
\[
\min(\sup_{s_1\in S_1}f(s_1),\sup_{s_2\in S_2}f(s_2))+\max(\sup_{s_1\in S_1}f(s_1),\sup_{s_2\in S_2}f(s_2))=\sup_{s_1\in S_1}f(s_1)+\sup_{s_2\in S_2}f(s_2)
\]
yields after integration with respect to $\sigma(df)$
\[
\int_{\bbsC_0}\min(\sup_{s_1\in S_1}f(s_1),\sup_{s_2\in S_2}f(s_2))\,\sigma(df)+\theta(S_1\cup S_2)=\theta(S_1)+\theta(S_2).
\]
This proves the first point of Theorem  \ref{theo2}. 
\item The second point is straightforward since
\begin{eqnarray*}
\bbP[f\not<_{S_1}\eta,\ f\not<_{S_2}\eta]&=&\bbP[\exists i\in I,\ f\not<_{S_{1,i}}\eta,\ \exists j\in J,\ f\not<_{S_{2,j}}\eta]\\
&\leq&\sum_{i\in I}\sum_{j\in J}\bbP[f\not<_{S_{1,i}}\eta,\ f\not<_{S_{2,j}}\eta]
\end{eqnarray*}
so that
\begin{eqnarray*}
\beta(S_1,S_2)&\leq& \sum_{i\in I}\sum_{j\in J}\int_{\bbsC_0}\bbP[f\not<_{S_{1,i}}\eta,\ f\not<_{S_{2,j}}\eta]\, \mu(df)\\
&\leq & 2\sum_{i\in I}\sum_{j\in J}\big[C(S_{1,i})+C(S_{2,j})\big]\ \big[\theta(S_{1,i})+\theta(S_{2,j})-\theta(S_{1,i}\cup S_{2,j})\big].
\end{eqnarray*}
\end{itemize}
\end{preuve}

\begin{preuve}{of Corollary \ref{cor2}}
This is a straightforward consequence of the second point of Theorem \ref{theo2} with $S_1=\cup_{s_1\in S_1}\{s_1\}$ and $S_2=\cup_{s_2\in S_2}\{s_2\}$. It holds indeed 
\[
\beta(S_1,S_2)\leq 2\sum_{s_1\in S_1}\sum_{s_2\in S_2}\big[C(\{s_{1}\})+C(\{s_{2}\})\big]\ \big[\theta(\{s_{1}\})+\theta(\{s_{2}\})-\theta(\{s_{1}\}\cup \{s_{2}\})\big]
\]
with 
\[
\theta(\{s_{1}\})=\theta(\{s_{2}\})=1\quad ,  \quad \theta(\{s_{1}\}\cup \{s_{2}\})=\theta(s_1,s_2)
\]
and
\[
 C(\{s_{1}\})=C(\{s_{2}\})=1.
\]
The last equality follows from the fact that, for all $s\in S$, $\eta(s)$ has a standard unit Fr\'echet distribution and hence $\eta(s)^{-1}$ has an exponential distribution with mean $1$. Hence we obtain
\[
\beta(S_1,S_2)\leq 4\sum_{s_1\in S_1}\sum_{s_2\in S_2}(2-\theta(s_{1},s_{2}))
\]

\end{preuve}

\subsection{Proof of Theorems \ref{theo3} and Proposition \ref{prop:asynorm}}
\begin{preuve}{of Theorem \ref{theo3}}
According to Bolthausen's CLT for stationary mixing random fields (see Appendix~\ref{sec:TCL}), it is enough to prove that the mixing coefficients $\alpha_{k,l}(m)$ defined by Equation \eqref{eq:amix} with $X(t)=g(\eta(t_1+t),\ldots,\eta(t_p+t))$ satisfy Equations \eqref{eq:TCL1}, \eqref{eq:TCL2} and \eqref{eq:TCL3}.

For $S\subset\bbZ^d$, we define $\widetilde S=\cup_{i=1}^p\{s+t_i,\ s\in S\}$. The inclusion of $\sigma$-fields
\[
\sigma(\{X(s),\ s\in S\})\subset \sigma(\{\eta(s),\ s\in \widetilde S\})
\]
entails a comparison of the related $\alpha$-mixing coefficients: for disjoint $S_1,S_2\subset \bbZ$,
\[
\alpha^{X}(S_1,S_2)\leq \alpha^{\eta}(\widetilde S_1,\widetilde S_2),
\]
where the superscript $X$ or $\eta$ denotes that we are computing the $\alpha$-mixing coefficient of the random field $X$ or $\eta$ respectively. Furthermore, 
\[
|\widetilde S_i|\leq p|S_i|, \ i=1,2\quad \mbox{and}\quad d(\widetilde S_1,\widetilde S_2)\geq d(S_1,S_2)-\Delta,
\]
with $\Delta=\max_{1\leq i<j\leq p} d(t_i,t_j)$  the diameter of $\{t_1,\dots,t_p\}$.
Hence, with obvious notations,
\[
\alpha_{k,l}^X(m)\leq \alpha_{pk,pl}^\eta(m-\Delta),\quad k,l\in\bbN\cup\{\infty\},\ m\geq 1.
\]
Corollary \ref{cor2} implies 
\[
\alpha_{k,l}^X(m)\leq \alpha_{pk,pl}^\eta(m-|h|)\leq p^2kl \sup_{|t|\geq m-\Delta}\gamma(t),\quad  k,l\in\bbN,\ m\geq 1,
\]
and 
\[
\alpha_{1,\infty}^X(m)\leq \alpha_{p,\infty}^\eta(m-\Delta)\leq p \sum_{|t|\geq m-\Delta}\gamma(t),\quad m\geq 1.
\]
In view of this, Assumption \eqref{eq:condgamma} entails Equations \eqref{eq:TCL1}, \eqref{eq:TCL2} and \eqref{eq:TCL3}, so that the random field $X$ satisfies Bolthausen's CLT.
\end{preuve}

\begin{preuve}{of Proposition \ref{prop:asynorm}}
Let $h\in\bbZ^d$. We apply Theorem \ref{theo3} to the stationary random field 
\[
X(t)=1_{\{\eta(t)\leq y,\ \eta(t+h)\leq y\}},\quad t\in\bbZ^d.
\]
Clearly $\bbE[X(t)]=\exp[-\theta(h)/y]$ and  $\bbE[|X|^{2+\delta}]<\infty$ for all $\delta>0$.
Assumption \eqref{eq:asynorm1} together with $\gamma(t)\leq 4(2-\theta(t))$ ensure that Equation \eqref{eq:condgamma} is satisfied for $\delta$ large enough. Hence the estimator $\hat p_n(h,y)$ is asymptotically normal:
\[
|\Lambda_n|^{1/2}\Big(\hat p_n(h,y)-p(h,y)\Big)\Longrightarrow \cN(0,\beta_1^2)
\]
with limit  variance 
\begin{eqnarray*}
\beta_1^2(y)&=&\sum_{t\in\bbsZ^d}\mathrm{Cov}[X(0),X(t)]\\
&=&\sum_{t\in\bbsZ^d}\Big( \exp(\theta(\{0,h,t,t+h\})/y)-\exp(2\theta(h))/y)\Big)\ >0.
\end{eqnarray*}
The $\delta$-method  entails  the asymptotic normality of the estimator $\hat\theta_n^y(h)=-y\log\hat p_n(h,y)$:
\[
|\Lambda_n|^{1/2}\Big(\hat \theta_n^y(h)-\theta(h)\Big)\Longrightarrow \cN(0,\sigma^2_1)
\]
with limit variance
\begin{eqnarray*}
\sigma^2_1=y^{2}\exp(2\theta(h)/y)\beta^2_1=y^2\sum_{t\in\bbsZ^d}\Big(\exp[(2\theta(h)-\theta(\{0,h,t,t+h\}))/y]-1\Big).
\end{eqnarray*}
The proof of the asymptotic normality of $\hat\theta_n^{(2)}$ and $\hat\theta_n^{(3)}$ is very similar and we give only the main lines. Using Theorem \ref{theo3}, we prove that
\[
\hat q_n= |\Lambda_n|^{-1}\sum_{t\in|\Lambda_n|} \min(\eta(t)^{-1},\eta(t+h)^{-1})
\]
is an asymptotic normal estimator of $\theta(h)^{-1}$:
\[
|\Lambda_n|^{1/2}(\hat q_n-\hat\theta(h)^{-1})\Longrightarrow \cN(0,\beta_2^2)
\]
with limit variance
\[
\beta_2^2=\sum_{t\in\bbsZ^d}\mbox{Cov}\big[\min(\eta(0)^{-1},\eta(h)^{-1}),\min(\eta(t)^{-1},\eta(t+h)^{-1})\big].
\]
The $\delta$-method entails the asymptotic normality of $\theta_n^{(2)}(h)=1/\hat q_n(h)$ with limit variance
\[
\sigma_2^2=\theta(h)^4\beta_1^2.
\]
Similarly, 
\[
\hat\nu_{F,n}(h)=|\Lambda_n|^{-1}\sum_{t\in|\Lambda_n|} |F(\eta(t))-F(\eta(t+h))|
\]
is an asymptotic normal estimator of $\nu_F(h)=\bbE[|F(\eta(0))-F(\eta(h))|]$:
\[
|\Lambda_n|^{1/2}(\hat \nu_{F,n}(h)-\nu_F(h))\Longrightarrow \cN(0,\beta_3^2)
\]
with limit variance
\[
\beta_3^2=\sum_{t\in\bbsZ^d}\mbox{Cov}\big[|F(\eta(0))-F(\eta(h))|,|F(\eta(t))-F(\eta(t+h))|\big].
\]
The $\delta$-method entails the asymptotic normality of 
\[
\theta_n^{(3)}(h)=\frac{1+2\hat\nu_{F,n}(h)}{1-2\hat\nu_{F,n}(h)}
\]
 with limit variance
\[
\sigma_3^2=(\theta(h)+1)^4\beta_3^2.
\]
\end{preuve}

\appendix
\section{Auxiliary results}
\subsection{Structure of max-i.d. random processes}\label{sec:A1}
The structure of sample continuous random processes on a compact metric space  was elucidated by Hahn, Gin\'e and Vatan \cite{GHV90}. We extend here their results for $T$  a locally compact metric space, in order to cover the standard cases $T=\bbZ^d$ or $T=\bbR^d$. Such extensions have been considered for max-stable models on $\bbR$ (see \cite{dHF06} Chapter 9.6) but we have found no reference in the max-i.d. case. 

Let $\eta$ be a continuous max-i.d. random process on $\bbC(T,\bbR)$. Define its vertex function $h:T\to [-\infty,+\infty)$ by
\[
h(t)=\mbox{essinf}\ \eta(t)=\sup\{x\in\bbR;\ \bbP(\eta(t)\geq x)=1\}.
\]
We will always assume that $h$ is continuous. We can then suppose without loss of generality that $h\equiv 0$. Indeed,  if $h$ is continuous and finite, we may consider $\eta-h$ which is a continuous max-i.d. random field with zero vertex function; and if $h$ is not finite everywhere, we may consider $\exp(\eta)-\exp(h)$ which is max-i.d with zero vertex function.

We note $\bbC(T)=\bbC(T,[0,+\infty))$ the space of nonnegative continuous function on $T$ endowed with the topology of uniform convergence on compact sets and $\bbC_0(T)=\bbC(T)\setminus\{0\}$.
\begin{theo}\label{theo:maxid} 
\begin{itemize}
\item Let $\eta=(\eta(t))_{t\in T}$ be a continuous max-i.d. process  on $T$ with vertex function $h\equiv 0$. There exists a unique locally finite Borel measure on $\bbC_0$ satisfying condition \eqref{eq:condmu}, called the exponent measure of $\eta$, such that
\[
\mu\Big[\cup_{i=1}^k\{f\in\bbC_0; f(t_i)>y_i\} \Big]=-\log\bbP\Big[\cap_{i=1}^k\{\eta(t_i)\leq y_i\}\Big]
\]
for all $k\geq 1$, $t_1,\ldots,t_k\in T^k$ and $y_1,\ldots,y_k>0$. 
\item Conversely, for any locally  finite Borel measure on $\bbC_0$ satisfying condition \eqref{eq:condmu}, there exists a continuous max-i.d. process $\eta$ on $T$ with vertex function $h\equiv 0$ and exponent measure $\mu$. It can be constructed as follows: let $\Phi$ be a Poisson point process on $\bbC_0$ with intensity $\mu$ and define 
\[
\eta(t)=\max\{\phi(t),\ \phi\in\Phi\},\quad t\in T.
\]
\end{itemize}
\end{theo}
\begin{preuve}{of Theorem \ref{theo:maxid}}
Let $(T_n)_{n\geq 1}$ be an increasing sequence of compact sets such that $T=\cup_{n\geq 1} T_n$. We suppose furthermore that $T_{n}$ is included in the interior set of $T_{n+1}$. The space $\bbC(T)=\bbC(T,\bbR^+)$ of nonnegative continuous functions on $T$ endowed with the topology of uniform convergence on compact sets can be seen as the projective limit of the sequence of spaces $\bbC(T_n,\bbR^+)$ endowed with the topology of uniform convergence. 
For $m\geq n\geq 1$, we define the natural projections 
\[
\pi_n:\bbC(T)\to\bbC(T_n)\quad \mbox{and}\quad \pi_{n,m}:\bbC(T_m)\to\bbC(T_n).
\]
For each $n\geq 1$, the restriction $\pi_n(\eta)=\eta_{T_n}$ is a continuous max-i.d. process on the compact space $T_n$ and according to \cite{GHV90}, there exists a locally finite exponent measure $\mu_n$ on $\bbC_0(T_n)=\bbC(T_n)\setminus\{0\}$ satisfying equation 
\[
\mu_n\Big[\cup_{i=1}^k\{f\in\bbC_0; f(t_i)>y_i\} \Big]=-\log\bbP\Big[\cap_{i=1}^k\{\eta(t_i)\leq y_i\}\Big]
\]
for all $k\geq 1$, $t_1,\ldots,t_k\in T_n^k$ and $y_1,\ldots,y_k>0$. 
Furthermore, for all $\varepsilon>0$
\[
\mu_n[\cS_{n,\varepsilon}]<\infty\quad \mbox{where}\quad \cS_{n,\varepsilon}=\Big\{f\in\bbC(T_n);\ \sup_{T_n} f> \varepsilon\Big\}.
\]
Let $n_0\geq 1$ and $\varepsilon>0$ be fixed. For $n\geq n_0$, define the  finite Radon measure by
\[
\tilde\mu_n^{n_0,\varepsilon}[A]=\mu_n[A\cap \pi_{n_0,n}^{-1}\cS_{n_0,\varepsilon}],\quad A\subset \bbC(T_n)\ \mbox{Borel\ set}.
\]
Clearly, the following compatibility conditions holds true: for $m\geq n\geq n_0$,
\[
\tilde\mu_n^{n_0,\varepsilon}=\tilde\mu_m^{n_0,\varepsilon}\pi_{n,m}^{-1}. 
\]
Theorem 5.1.1 in \cite{B55} state the existence of projective limit of Radon measures, it implies the existence of a finite Radon measure $\tilde \mu^{n_0,\varepsilon}$ on $\bbC(T)$ such that
\[
\tilde\mu_n^{n_0,\varepsilon}=\tilde\mu^{n_0,\varepsilon}\pi_{n}^{-1},\quad n\geq n_0. 
\]
It is then easily checked  that the measure  $\mu$  on $\bbC_0(T)$ defined by 
\[
\mu[A]=\sup\{ \tilde\mu^{n_0,\varepsilon}[A];\ n_0\geq 1,\ \varepsilon>0\},\quad A\subset \bbC_0(T)\ \mbox{Borel\ set}
\]
is locally finite  and  enjoys the required properties.
\end{preuve}

\subsection{Proof of Lemma \ref{lem:lc}}\label{sec:lemlc}
First, please note that the event $\{f<_S\eta\}$  depends only on the restriction $\eta_{S}$ and is hence measurable with respect to  the $\sigma$-field generated by $\Phi_S^+$.\\
In order to prove the Proposition, let $A_1,\ldots,A_k\subset\bbC_0$ be disjoint compact sets and $n_1,\ldots,n_k\geq 0$. Let $A=\cup_{i=1}^k A_i$ and $n=\sum_{i=1}^k n_i$. We compute the conditional probability with respect to $\Phi_S^+$ of the event 
\[
\{\Phi_S^-(A_1)=n_1,\ldots,\Phi_S^-(A_k)=n_k\}.
\]
This event is equal to $\{\Phi_S^-\in B\}$ with $B=\{N\in M_p(\bbC_0);\ N(A_1)=n_1,\ldots,N(A_k)=n_k\}$.
We remark that that it is realized if and only if there exists a $n$-uplet $(\phi_1,\ldots,\phi_n)$ of  atoms of $\Phi$ such that:
\begin{itemize}
\item[-] the atoms $\phi_1,\ldots,\phi_n$ are $S$-subextremal;
\item[-] $\sum_{j=1}^n \delta_{\phi_j}\in B$;
\item[-] the point measure $\Phi-\sum_{j=1}^n \delta_{\phi_j}$  has no $S$-subextremal atom in $A$, i.e. it belongs to
\[
D=\{N\in M_p(\bbC_0);\ N_S^-(A)=0\}. 
\]
\end{itemize}
Then the $n$-uplet $(\phi_1,\ldots,\phi_n)$ is unique up to permutation of the coordinates. The above observations entail that for all measurable $C\subset M_p(\bbC_0)$,
\begin{eqnarray*}
\bbP[\Phi_{S}^+\in C,\ \Phi_S^-\in B]
&=&\frac{1}{n!}\bbE\Big[\int_{\bbsC_0^n} 1_{\{\Phi_S^+\in C\}}1_{\{\forall i\in[\![1,n]\!],\ \phi_i<_S \eta\}}1_{\{\sum_{i=1}^n\delta_{\phi_i}\in B\}}1_{\{\Phi-\sum_{i=1}^n\delta_{\phi_i}\in D\}}\quad\\
&&\qquad\qquad \Phi(d\phi_1)(\Phi-\delta_{\phi_1})(d\phi_2)\cdots(\Phi-\sum_{i=1}^{n-1}\delta_{\phi_i})(d\phi_n)\Big].
\end{eqnarray*}
Slyvniak's formula entails
\begin{eqnarray}
&&\bbP[\Phi_{S}^+\in C,\ \Phi_S^-\in B]\nonumber\\
&=&\bbE\Big[1_{\{\Phi_S^+\in C\}}1_{\{ \Phi_S^-(A)=0\}}\frac{1}{n!}\int_{\bbsC_0^n} 1_{\{\sum_{i=1}^n\delta_{f_i}\in B\}}\otimes_{i=1}^n \big(1_{\{f_i<_S\eta\}}\mu(df_i)\big)\Big]\label{eq:Slyvniak}.
\end{eqnarray}
Summing this relation over the different values of $n_1,\ldots,n_k\in\bbN$ and the related sets $B=\{N\in M_p(\bbC_0);\ N(A_1)=n_1,\ldots,N(A_k)=n_k\}$, we obtain
\[
\bbP[\Phi_{S}^+\in C]=\bbE\Big[1_{\{\Phi_S^+\in C\}}1_{\{ \Phi_S^-(A)=0\}}\exp\big[\mu(\{f\in A;\  f<_S\eta\})\big]\Big].
\]
So we can rewrite Equation \eqref{eq:Slyvniak} as
\begin{eqnarray*}
\bbP[\Phi_{S}^+\in C,\ \Phi_S^-\in B]&=&\bbE\Big[1_{\{\Phi_S^+\in C\}}1_{\{ \Phi_S^-(A)=0\}}\exp\big[\mu(\{f\in A;\ f<_S\eta\})\big]K(\eta_{S},B)\Big],
\end{eqnarray*}
where 
\[
K(\eta_{S},B)=\frac{\exp\big[-\mu(\{f\in A;\ f<_S\eta\})\big]}{n!}\int_{A^n} 1_{\{\sum_{i=1}^n\delta_{f_i}\in B\}}\otimes_{i=1}^n \big(1_{\{f_i<_S\eta\}}\mu(df_i)\big)
\]
 is the conditional probability of $\{\Phi_S^{-}\in B\}$ with respect to $\Phi_S^+$ (note it depends on $\Phi_S^+$ only through the restriction $\eta_{S}$). We recognize the distribution of a Poisson random measure with intensity $1_{\{f<_S\eta\}}\mu(df)$ and this proves Lemma \ref{lem:lc}.

\subsection{Slyvniak's formula}\label{sec:slyvniak}
Palm Theory deals with conditional distribution for point
processes. We recall here one of the most famous formula of Palm
theory, known as Slyvniak's Theorem. This will be the main tool in our
computations. For a general reference on Poisson point processes, Palm
theory and their applications, the reader is invited to refer to the
monograph \cite{L13} by Stoyan, Kendall and Mecke.

Let $M_{p}(\bbC_0)$ be the set of locally-finite point measures $N$ on $\bbC_0$ endowed with the $\sigma$-algebra generated by the family of mappings $\{N\mapsto N(A),\  A\subset \bbC_0 \mathrm{\ Borel\ set}\}$. 

\begin{theo}[Slyvniak's Formula]\ \\
Let $\Phi$ be a Poisson point process on $\bbC_0$ with intensity measure $\mu$. For all measurable function $F:\bbC_0^k\times M_{p}(\bbC_0)\to [0,+\infty)$, 
\begin{eqnarray*}
& &\bbE\Big[\int_{\bbsC_0^k} 
  F\Big(\phi_1,\ldots,\phi_k,\Phi-\sum_{i=1}^k \delta_{\phi_i}\Big)\,\Phi(d\phi_1)\,(\Phi-\delta_{\phi_1})(d\phi_2)\cdots (\Phi-\sum_{j=1}^{k-1} \delta_{\phi_j})(d\phi_k) \Big]\\
 &=&\int_{\bbsC_0^k}\bbE[F(f_1,\ldots,f_k,\Phi)]\,\mu^{\otimes k}(df_1,\ldots,df_k).
\end{eqnarray*}
\end{theo}

\subsection{A central limit Theorem for weakly dependent process}\label{sec:TCL}
Since the pioneer work of Ibragimov \cite{I62}, many versions of the central limit Theorem for weakly dependent processes have been developed under various strong mixing conditions. We present here a central limit Theorem for stationary mixing random fields due to Bolthausen \cite{B82}. Let $(X_k)_{k\in\bbsZ^d}$ be a real valued stationary random field and recall the definition of the $\alpha$-mixing coefficient \eqref{eq:amix}. 
If $\Lambda\subset\bbZ^d$, we note $|\Lambda|$ the number of elements in $\Lambda$ and $\partial\Lambda$ the set of elements $k\in\Lambda$ such that there is $l\notin \Lambda$ with $d(k,l)=1$. Let $\Lambda_n$ be a fixed increasing sequence of finite subsets of $\bbZ^d$, which increases to $\bbZ^d$ and such that $\lim_{n\to\infty}|\partial \Lambda_n|/|\Lambda_n|=0$.
Let $S_n=\sum_{h\in\Lambda_n} (X_h-\bbE[X_h])$.

Bolthausen's central limit Theorem is based on the mixing coefficients
\begin{equation}\label{eq:amix}
\alpha_{kl}(m)=\sup\Big\{\alpha(S_1,S_2);\ |S_1|=k, |S_2|=l,\ d(S_1,S_2)\geq m \Big\}
\end{equation}
defined for $m\geq 1$ and $k,l\in\bbN\cup\{\infty\}$.
\begin{theo}\label{theo:BCLT}
Suppose that the following three conditions are satisfied:
\begin{eqnarray}
&\alpha_{1\infty}(m)=o(m^{-d});\label{eq:TCL1}\\
&\sum_{m=1}^\infty m^{d-1}\alpha_{kl}(m)<\infty\quad  \mbox{for\ all\ \ } k\geq 1,\ l\geq 1\ \mbox{such\ that\ } k+l\leq 4;\label{eq:TCL2}\\
&\bbE\big[|X_h|^{2+\delta}\big]<\infty \quad\mbox{and}\quad \sum_{m=1}^\infty m^{d-1}|\alpha_{11}(m)|^{\delta/(2+\delta)}<\infty\quad\mbox{for\ some}\ \delta>0.\label{eq:TCL3}
\end{eqnarray}
Then the series $\sigma^2=\sum_{h\in\bbsZ^d} \mathrm{Cov}[X_0,X_h]$ converges absolutely and if furthermore $\sigma^2>0$,
\[
 \frac{S_n}{\sigma|\Lambda_n|^{1/2}}\Longrightarrow \cN(0,1),\quad \mbox{as}\ n\to\infty.
\]
\end{theo}

\bibliographystyle{plain}
\bibliography{Biblio2}

\end{document}